\documentclass[a4paper,12pt]{amsart}

\usepackage{amssymb,amscd,amsthm,mathrsfs}
\usepackage{epsfig,graphicx}

\usepackage[ansinew]{inputenc}

\usepackage[centertags]{amsmath}

   \def\DD{{\mathbb D}}

 \def\RR{{\mathbb R}} \def\SS{{\mathbb S}} \def\TT{{\mathbb T}}
   
 \def\ZZ{{\mathbb Z}}

\def\cA{\mathcal{A}}  \def\cG{\mathcal{G}}  
   \def\cN{\mathcal{N}} 
\def\cC{\mathcal{C}}   \def\cO{\mathcal{O}} \def\cU{\mathcal{U}}
    \def\cV{\mathcal{V}}
\def\cE{\mathcal{E}}    
\def\cF{\mathcal{F}}

\newcommand{\en}{\subset}

\newcommand{\T}{\mbox{$\mathbb{T}$}}
\newtheorem*{teorema}{Theorem}

\newtheorem*{teoA}{Theorem A}
\newtheorem*{teoB}{Theorem B}
\newtheorem{teo}{Theorem}[section]

\newtheorem{quest}{Question}
\newtheorem{cor}[teo]{Corollary}

\newtheorem{lema}[teo]{Lemma}
\newtheorem{prop}[teo]{Proposition}

\newcommand{\bi}{\begin{itemize}}
\newcommand{\ei}{\end{itemize}}

\theoremstyle{definition}

\theoremstyle{remark}
\newtheorem{obs}[]{Remark}

\newcommand{\Fol}{\mbox{$\mathcal{F}$}}
\newcommand{\D}{\mbox{$\mathbb{D}$}}

\newcommand{\eps}{\varepsilon}

\newcommand{\esbozo}{\vspace{.05in}{\sc\noindent Sketch. }}
\newcommand{\dem}[1]{\vspace{.05in}{\sc\noindent Proof #1.}}
\newcommand{\lqqd}{\par\hfill {$\Box$} \vspace*{.05in}}
\newcommand{\finobs}{\par\hfill{$\diamondsuit$} \vspace*{.05in}}

\newcommand{\U}{\mathcal{U}}
\newcommand{\trans}{\mbox{$\,{ \top} \;\!\!\!\!\!\!\raisebox{-.3ex}{$\cap$}\,$}}

\DeclareMathOperator{\Diff}{Diff}

\DeclareMathOperator{\diametro}{diam}

\author[Rafael Potrie]{Rafael Potrie}
\address{CMAT, Facultad de Ciencias, Universidad de la Rep\'ublica, Uruguay}
\address{LAGA; Institute Galilee, Universite Paris 13, Villetaneuse, France}
\email{rpotrie@cmat.edu.uy}

\title[Wild Milnor attractors accumulated by lower dimensional dynamics]{Wild Milnor attractors accumulated by lower dimensional dynamics}
\thanks{The autor was partially supported by ANR Blanc DynNonHyp BLAN08-2$\_$313375 and  ANII Proyecto FCE2007$\_$577}

\begin{document}

\maketitle


\begin{abstract}
We present new examples of open sets of diffeomorphisms such that
a generic diffeomorphisms in those sets have no dynamically
indecomposable attractors in the topological sense and have
infinitely many chain-recurrence classes. We show that except from
one particular class, the other classes are contained in periodic
surfaces. This study allows us to obtain existence of Milnor
attractors as well as studying ergodic properties of the
diffeomorphisms in those open sets by using the ideas and results
from \cite{BV} and \cite{BF}.
\end{abstract}

\section{Introduction}
\subsubsection{} In 1987, A. Araujo in his thesis (\cite{A}) announced that
$C^1$-generic diffeomorphisms of compact surfaces have hyperbolic
attractors. In fact, he claimed to have proved that for a residual
subset of diffeomorphisms on a compact surface, either there are
infinitely many sinks or there are finitely many hyperbolic
attractors whose basin cover a full Lebesgue measure of the
manifold. The proof seems to have a gap, but the techniques in
\cite{PS} allow to overcome them (and with the recent results of
$C^1$ generic dynamics this can be proven rather
easily\footnote{See \cite{Potrie}.}).

In contrast, an astonishing example was recently constructed in
\cite{BLY} where they showed that there exist open sets of
diffeomorphisms in any manifold of dimension $\geq 3$ such that
every $C^r$-generic diffeomorphism of those open subsets have no
attractors and there is an attracting region having infinitely
many distinct chain recurrence classes. We recommend reading the
introduction of \cite{BLY} for more on the history of this
important problem.

The construction in \cite{BLY} relies on some modification of the
well known solenoid attractor. Although the construction is rather
simple, it is not well understood how is that other chain
recurrence classes coexist in this attracting region.

\subsubsection{} In this paper we propose a new kind of example starting from a non
hyperbolic DA attractor (based on an example of \cite{Car}, see
also \cite{BV}) which allows us to use the properties of
semiconjugacy with a linear Anosov diffeomorphisms, and gives a
more satisfactory picture of how the infinitely many
chain-recurrence classes behave.

Also, this study allows us to obtain some remarkable features of
this examples from the ergodic point of view which can be
summarized in the following statement (to be stated in a precise
form in section \ref{SectionEnunciados}).

\begin{teorema}
There exists a $C^1$-open set $\cU$ of $\Diff^r(\TT^3)$ ($r\geq
1$) and a $C^r-$generic subset $\cG^r \en \cU$ such that every
$f\in \cG^r$ has no attractors and $f$ has infinitely many
chain-recurrence classes. Moreover:
\begin{itemize}
\item[-] For every $f\in \cU$ there exists a chain-recurrence
class $H$ such that every chain-recurrence class $R$ different
from $H$ is contained in a periodic surface. \item[-] For every
$f\in \cU$ there exists a unique attractor in the sense of Milnor.
\item[-] For every $f\in \cU$ there exists a unique entropy
maximizing measure. \item[-] If $r\geq 2$ and $f\in \cU$ then $f$
admits a unique SRB measure.
\end{itemize}
\end{teorema}

\subsubsection{} Besides this ergodic point of view, there is another motivation in
studying the examples here proposed. Recently, C. Bonatti has
proposed a program for studying $C^1$-generic dynamics (\cite{B})
and in particular, what is known as wild dynamics (see section
\ref{DefinicionBDP}). He has defined viral homoclinic classes as
those classes essentially having a reproductive behavior (until
now, the only known mechanism for generating wild dynamics, see
\cite{BD1, BCDG}). Even if we are not able to prove that the
examples here presented are not viral (we shall not define this
notion here, see \cite{B} or \cite{BCDG} for a precise
definition), the fact that all chain-recurrence classes except one
are contained in periodic surfaces seems to represent a different
mechanism for generating wild dynamics.

On the other hand, when studying wild homoclinic
classes\footnote{To be defined in section \ref{DefinicionBDP}.}
with a partially hyperbolic structure, it has been announced by C.
Bonatti and K. Shinohara that the examples from \cite{BLY} are
viral and it seems that the main feature differentiating the
behaviors is the topology of the intersections between the
homoclinic classes and the center-stable manifolds, this becomes
clear in our Proposition \ref{ProposicionMecanismo}.

\subsection{Some definitions and results which will be used}\label{SectionDefiniciones}

We shall give some definitions and state some results we shall use
along the paper, it may be wise to skip this section and return to
it when not knowing some definition or when it is referred to by
the text.

\subsubsection{}\label{DefinicionCONLEY} {\bf Conley's theory and Bonatti-Crovisier's result.}
Given a homeomorphism $f: M \to M$ we can define the following
relation on $M$: we denote $x \dashv y$ whenever for every
$\eps>0$ there exists an $\eps-$pseudo-orbit from $x$ to $y$, that
is, there exists a set of points $x=z_0, \ldots , z_n=y$ such that
$n\geq 1$ and $d(f(x_i),x_{i+1})<\eps$.

We denote as $$R(f)= \{ x \in M \ : \ x\dashv x \}$$ \noindent the
\emph{chain-recurrent set} of $f$. In $R(f)$ the relation $x
\vdash \!\!\! \dashv y$ (given by $x \vdash \!\!\! \dashv y$ if
and only if $x\dashv y$ and $y \dashv x$) is an equivalence
relation, we shall call its equivalence classes
\emph{chain-recurrence classes}. An invariant set will be called
\emph{chain-transitive} if it is transitive under the relation
$\vdash\!\!\! \dashv$.

An open set $U$ is a \emph{filtrating neighborhood} if there
exists $V_1, V_2$ open sets such that $U = V_1 \setminus \overline
V_2$ and $f(\overline{V_i}) \en V_i$ for $i=1,2$. It is not hard
to see that a filtrating neighborhood contains each
chain-recurrence class it intersects.

Conley's theorem (see \cite{Crov}) implies that given two
different chain-recurrence classes $\cC_1$ and $\cC_2$, there
exists a filtrating neighborhood containing $\cC_1$ which does not
intersect $\cC_2$. A chain-recurrence class $\cC$ is
\emph{isolated} if there exists a filtrating neighborhood $U$ such
that $U\cap R(f)= \cC$.

We shall pay special attention to certain particular
chain-recurrence classes: We say that a compact invariant set $Q$
is a \emph{quasi-attractor}  if it is a chain-recurrence class and
there exists a decreasing sequence of open neighborhoods $\{U_n\}$
such that $\bigcap U_n = Q$ and $f(\overline{U_n})\en U_n$. An
important feature is that for quasi-attractors one has that if a
point $y$ verifies that there exists $x\in Q$ such that $x\dashv
y$, then $y \in Q$, in particular, $Q$ is saturated by unstable
sets.

We say that a compact invariant set $Q$ is an
\emph{attractor}\footnote{It is more usual to find in the
literature the following definition: A compact $f$-invariant set
$\Lambda$ is an \emph{attractor} if it is contains a dense orbit
and there is a neigbhorhood $U$ of $\Lambda$ such that
$f(\overline U)\en U$ and $\Lambda= \bigcup_{n\geq 0}
f^n(\overline U)$. Our definition coincides with this one except
that we demand the weaker indecomposability hypothesis of being
chain-transitive instead of transitive.} if it is a
quasi-attractor and is isolated as chain-recurrence class.

For a diffeomorphism $f$ of $M$ and a point $x\in M$ we define
$$W^s(x)= \{ y\in M \ : \ d(f^n(x),f^n(y)) \to 0 \, \ n\to +\infty
\}$$ \noindent and $$W^u(x) = \{ y\in M \ : \ d(f^n(x),f^n(y)) \to
0 \, \ n\to -\infty \}$$ \noindent the stable and unstable sets of
$x$.

For hyperbolic periodic points, it is well known that this sets
are $C^1$ injectively immersed manifolds, the \emph{stable index}
of a periodic point is the dimension of its stable manifold. We
define the \emph{homoclinic class of a periodic point $p$} as the
closure of the transversal intersections between the stable and
unstable set of the points in the orbit of $p$ (i.e. $H(p)=
\overline{W^s(\cO(p)) \trans W^u(\cO(p))}$).

If for a point $x\in M$ we have that $W^\sigma (x)$ is a manifold,
we shall denote as $W^\sigma_L (x)$ as the disk of radius $L$
centered at $x$ in $W^\sigma(x)$ with the Riemannian metric
induced in $W^\sigma(x)$ by the immersion.

It was proved in \cite{BC} that for a residual ($G_\delta$-dense)
subset $\cG_{BC}$ of $\Diff^1(M)$ one has the following
properties:

\begin{itemize}
\item[-] Every periodic point of $f$ is hyperbolic (this is the well known Kupka-Smale's theorem).
\item[-] $R(f) = \overline{Per(f)}$ where $Per(f)$ denotes the set of periodic points of $f$.
\item[-] If a chain-recurrence class $\cC$ of $f$ contains a periodic point $p$, then $\cC$ coincides with its homoclinic class $H(p)$.
\item[-] For a residual set of points $G \en M$ the omega-limit set is a quasi-attractor.
\end{itemize}

Since homoclinic classes are always transitive, we get from this
result that for $C^1$-generic diffeomorphisms (this will stand for
diffeomorphisms in a residual subset of $\Diff^1(M)$) the
definition of attractor we gave coincides with the usual one.

\subsubsection{}\label{DefinicionPartialHyp} {\bf Partial hyperbolicity.}
Given a diffeomorphism $f : M\to M$, a compact $f$-invariant
subset $\Lambda \en M$ and two $Df$-invariant subbundles $E$ and
$F$ of $T_\Lambda M$ we say that $F$ \emph{dominates} $E$ if there
exists $N>0$ such that for every $x\in \Lambda$ and every pair of
unit vectors $v\in E$, $w\in F$ we have that:

$$ \|Df_x^N v \| < \frac 1 2 \|Df_x^N w \| $$

Whenever there exists a cone field $\cE$ (of constant dimension)
such that for every $x\in \Lambda$ one has that $D_xf \cE(x) \en
int(\cE(f(x)))$ there exists a (unique) $Df$-invariant sum
$T_\Lambda M= E\oplus F$ with $F \en \cE$ and $E \en \cE^c$ and
such that $F$ dominates $E$ (see \cite{BDV} appendix B).

We say that $\Lambda$ is \emph{partially hyperbolic} provided that
$T_\Lambda M$ decomposes as a $Df$-invariant sum $T_\Lambda M =
E^{cs} \oplus E^u$ and there exists $N>0$ such that the following
conditions hold:

\begin{itemize}
\item[-] $E^u$ dominates $E^{cs}$.
\item[-] For every $x\in \Lambda$ and every unit vector $v^u$ in $E^u(x)$ we have that $\|Df_x^N v^u \| > 2$.
\end{itemize}

The definitions found in the literature (see particularly
\cite{BDV} appendix B or \cite{Crov}) require that either $f$ or
$f^{-1}$ is partially hyperbolic under our definition. We shall
not be concerned with this fact since it will be clear how to
adapt the results here to that definition.

In general (see \cite{BDV} appendix B or \cite{HPS}) we have that
the bundle $E^u$ integrates into a $f-$invariant lamination
$\cF^u$ of leaves tangent to $E^u$ which we shall call
\emph{strong unstable manifolds}.

When $E^{cs}$ also integrates into a $f$-invariant lamination
$\cF^{cs}$ tangent to $E^{cs}$ at $\Lambda$ we shall say that
$E^{cs}$ is \emph{coherent}.

The leaf of $\cF^{\sigma}$ through $x$ shall be denoted as
$\cF^\sigma(x)$ and $\cF^\sigma_L(x)$ will denote the ball of
radius $L$ centered at $x$ in $\cF^\sigma(x)$ with the induced
metric ($\sigma= u, cs$). Notice that for every point  $x\in
\Lambda$ we have that $\cF^u(x) \en W^u(x)$.

By \emph{lamination} on a set $K$ we mean a collection of disjoint
$C^1$ injectively immersed manifolds of the same dimension (called
\emph{leaves}) such that there exists a compact metric space
$\Gamma$ such that for every point $x\in K$ there exists a
neighborhood $U$ and a homeomorphism $\varphi: U \cap K \to \Gamma
\times \RR^d$ such that if $L$ is a leaf of the lamination and
$\tilde L$ a connected component of $L\cap U$ then
$\varphi|_{\tilde L}$ is a $C^1$-diffeomorphism to $\{s\} \times
\RR^d$ for some $s \in \Gamma$.

When the laminated set $K$ is the whole manifold, we say that the lamination is a \emph{foliation}.

Consider a compact invariant set $\Lambda$ which is partially
hyperbolic with splitting $E^{cs}\oplus E^u$ and such that $f$ is
coherent in $\Lambda$. Given an open set $U$ of $\Lambda$ and
$x\in \Lambda$, we denote as $\cF^{cs}_U(x)$ to the connected
component of $\cF^{cs}(x)\cap U$ containing $x$. If $y \in
\cF^u(x)$ is close to $x$ we can define the \emph{unstable
holonomy} $\Pi^{uu}_{x,y}$ from a neighborhood of $x$ in
$\cF^{cs}_U(x) \cap \Lambda$ to a neighborhood of $y$ in
$\cF^{cs}_U(y) \cap \Lambda$ as projecting the points along the
unstable leaves, which is a continuous injective map.

\subsubsection{}\label{DefinicionBDP} {\bf Non-isolated classes for the $C^1$-topology.}
We say that a diffeomorphism $f \in \Diff^1(M)$ is \emph{tame} if
it belongs to the interior of the diffeomorphisms having only
finitely many chain-recurrence classes. A diffeomorphism is
\emph{wild} if it cannot be accumulated by tame diffeomorphisms in
the $C^1$-topology. We get that in $\Diff^1(M)$ the union of tame
and wild diffeomorphisms is open and dense.

When $f$ is wild and $C^1$-generic, as a consequence of the
results of \cite{BC} we have that it has infinitely many
chain-recurrence classes, in particular, it will have at least one
which is not isolated. When a homoclinic class $H(p)$ is not
isolated, we shall say that it is a \emph{wild homoclinic class}.

As a consequence of the main result in \cite{BDP}, we obtain the
following criterium for partially hyperbolic homoclinic classes to
be wild:

\begin{teo}[\cite{BDP}]\label{BDP}
There exists a residual subset $\cG_{BDP}$ of $\Diff^1(M)$ such
that if a homoclinic class $H(p)$ verifies that:
   \begin{itemize}
   \item[-] The homoclinic class $H(p)$ admits a partially hyperbolic splitting of the form $T_{H(p)}M = E^{cs} \oplus E^u$.
   \item[-] The subbundle $E^{cs}$ admits no decomposition into non-trivial $Df-$invariant subbundles which are dominated.
   \item[-] There is a periodic point $q\in H(p)$ such that $\det (Df^{\pi(q)}_q |_{E^{cs}(q)}) > 1$.
   \end{itemize}
Then, $H(p)$ is contained in the closure of the set of periodic
sources\footnote{Periodic repelling points. Notice that the
chain-recurrence class of a periodic source is reduced to the
source itself.} of $f$. In particular, $H(p)$ is a wild homoclinic
class.
   \end{teo}

\subsubsection{}\label{DefinicionPV} {\bf Robust tangencies and $C^r-$generic non-isolation.}
As well as in the case of $C^1$-topology, we can obtain a similar
criterium to obtain non-isolation of a homoclinic class for
$C^r-$generic diffeomorphisms combining the main results of
\cite{BD} and \cite{PV}. The only cost will be that we must
consider a new open set and that the accumulation by other classes
is not as well understood (in Theorem \ref{BDP} we obtained that
the class is contained in the closure of the sources, and here we
shall only obtain that the class intersects the closure of the
sources). We state a consequence of the results in those papers in
the following result. We shall only use the result in dimension 3,
so we state it in this dimension, it can be modified in order to
hold in higher dimension but it would imply defining sectionally
dissipative saddles (see \cite{PV}).

\begin{teo}[\cite{BD} and \cite{PV}]\label{BDPV}
Consider $f \in \Diff^r(M)$ with $\dim M =3$ and a $C^1-$open set
$\cU$ of $\Diff^r(M)$ ($r\geq 1$) such that there exists a
hyperbolic periodic point $p$ of $f$ such that its continuation
$p_g$ is well defined for every $g\in \cU$ and such that:
   \begin{itemize}
   \item[-] The homoclinic class $H(p_g)$ admits a partially hyperbolic splitting of the form $T_{H(p)}M = E^{cs} \oplus E^u$ for every $g\in \cU$.
   \item[-] The subbundle $E^{cs}$ admits no decomposition in non-trivial $Dg-$invariant subbundles which are dominated.
   \item[-] There is a periodic point $q\in H(p_f)$ such that $\det (Df^{\pi(q)}_q |_{E^{cs}(q)}) > 1$.
   \end{itemize}
Then, there exists a $C^1$-open and dense subset $\cU_1 \en \cU$
and a $C^r$-residual subset $\cG_{PV}$ of $\cU_1$ such that for
every $g\in \cG_{PV}$ one has that $H(p_g)$ intersects the closure
of the set of periodic sources of $g$.
\end{teo}

The conditions of the Theorem are used in \cite{BD} in order to
create robust tangencies for a hyperbolic set for diffeomorphisms
in an $C^1$-open and dense subset $\cU_1$ of $\cU$. Then, using
similar arguments as in \cite{BLY} (section 3.7) one creates
tangencies associated with periodic orbits which are sectionally
dissipative for $f^{-1}$ which allows to use the results in
\cite{PV} to get the conclusion.

\subsubsection{}\label{DefiMILNORSRB} {\bf Milnor attractors, SRB measures and entropy maximizing measures.}
Following \cite{Milnor}, we shall say that a compact invariant
chain transitive set $\Lambda$ is an \emph{attractor in the sense
of Milnor} or \emph{Milnor attractor} if and only if the basin
$B(\Lambda)$ of $\Lambda$ has positive Lebesgue measure and for
every $\tilde \Lambda \subsetneq \Lambda$ compact invariant set,
its basin $B(\tilde \Lambda)$ has strictly smaller measure.
Moreover, if for every $\tilde \Lambda \subsetneq \Lambda$ compact
invariant set, we have that $Leb(B(\tilde \Lambda))=0$, we will
say that $\Lambda$ is a \emph{minimal attractor in the sense of
Milnor} or \emph{minimal Milnor attractor}.

Here, \emph{basin} must be understood as the set of points whose
forward iterates converge to the compact set (it must not be
confused with the statistical basin which is quite more
restrictive).

We say that an invariant measure $\mu$ is an \emph{SRB-measure}
whenever its statistical basin has positive Lebesgue measure, this
means that there exists a positive Lebesgue measure set $B(\mu)$
such that for every $x\in B(\mu)$ one has that for every
continuous function $\varphi: M \to \RR$

    $$ \frac 1 n \sum_{i=0}^{n-1} \varphi (f^i(x)) \to \int \varphi d\mu $$

We refer the reader to \cite{BDV} chapter 11 (in particular
section 11.2) for a nice introduction to SRB measures in this
exact context. Notice that the existence of an ergodic SRB measure
implies the existence of a minimal Milnor attractor.

Finally, we recall that an ergodic invariant measure $\mu$ is
called \emph{entropy maximizing measure} for a homeomorphism $f$
whenever its measure-theoretical entropy coincides with the
topological entropy of $f$. In \cite{BF} (see also \cite{BFSV})
some conditions where studied which imply existence and uniqueness
of these measures, this conditions will be satisfied on our
example.

\subsection{Precise statement of results}\label{SectionEnunciados}
Now we are in conditions to state our main results.

\begin{teoA} There exists a $C^1$-open set $\cU$ of $\Diff^r(\TT^3)$ such that:
\begin{itemize}
\item[(a)] For every $f\in \cU$ we have that $f$ is partially
hyperbolic with splitting $TM= E^{cs} \oplus E^u$ and $E^{cs}$
integrates to a $f$-invariant foliation $\cF^{cs}$. \item[(b)]
Every $f\in \cU$ has a unique quasi-attractor $Q_f$ which contains
a homoclinic class. \item[(c)] Every chain recurrence class $R\neq
Q_f$ is contained in the orbit of a periodic disk in a leaf of the
foliation $\cF^{cs}$. \item[(d)] There exists a residual subset
$\cG^r$ of $\cU$ such that for every $f\in \cG^r$ the
diffeomorphism $f$ has no attractors. In particular, $f$ has
infinitely many chain-recurrence classes accumulating on $Q_f$.
\item[(e)] For every $f\in \cU$ there is a unique Milnor attractor
$\tilde Q \en Q_f$. \item[(f)] If $r\geq 2$ then every $f \in \cU$
has a unique SRB measure whose support coincides with a homoclinic
class. Consequently, $\tilde Q$ is a minimal attractor in the
sense of Milnor. If $r=1$, then there exists a residual subset
$\cG_M$ of $\cU$ such that for every $f\in \cG_M$ we have that
$\tilde Q$ coincides with $Q_f$ and is a minimal Milnor attractor.
\item[(g)] For every $f\in \cU$ there is a unique entropy
maximizing measure.
\end{itemize}
\end{teoA}

By inspection in the proofs, one can easily see that in fact the
construction can be made in higher dimensional torus, however, it
can only be done in the isotopy classes of Anosov diffeomorphisms.
Also, it can be seen that condition (d) can be slightly
strengthened in the $C^1$-topology (see section
\ref{SectionTeorema37}).

We are able to construct examples which can be embedded in every
isotopy class of diffeomorphisms of manifolds of dimension $\geq
3$ on which we were not able to obtain the same ergodic
properties:

\begin{teoB} For every $d$-dimensional manifold $M$ and every isotopy class of diffeomorphisms of $M$ there exists a $C^1$-open set $\cU$ of $\Diff^r(M)$ such that for some open neighborhood $U$ in $M$:
\begin{itemize}
\item[(a)] Every $f\in \cU$ has a unique quasi-attractor $Q_f$ in $U$ which contains a homoclinic class and has a partially hyperbolic splitting $T_{Q_f}M = E^{cs} \oplus E^u$ which is coherent.
\item[(b)] Every chain recurrence class $R\neq Q_f$ is contained in the orbit of a periodic leaf of the lamination $\cF^{cs}$ tangent to $E^{cs}$ at $Q_f$.
\item[(c)] There exists a residual subset $\cG^r$ of $\cU$ such that for every $f\in \cG^r$ the diffeomorphism $f$ has no attractors. In particular, $f$ has infinitely many chain-recurrence classes.
\item[(d)] For every $f\in \cU$ there is a unique Milnor attractor $\tilde Q \en Q_f$.
\end{itemize}
\end{teoB}

The examples here are modifications of the product of a Plykin
attractor and the identity on the circle $\SS^1$. One can also
obtain them in order to provide examples of robustly transitive
attractors in dimension $3$ with splitting $E^{cs} \oplus E^u$
(see Appendix A). The author is not aware of other known examples
of such attractors other than Carvalho's example which is only
possible to be made in certain isotopy classes of diffeomorphisms.

\subsection{Further remarks on the construction and some questions}
\subsubsection{} As we mentioned, for a $C^1$-generic wild diffeomorphism, there
are infinitely many chain-recurrence classes. However, it is not
known if the cardinal of classes must be necessarily uncountable
(this holds in all examples where we know the cardinal of the
classes).

We pose the following question:

\begin{quest} In the open set of Theorem A, does it hold that for $C^1$-generic diffeomorphisms in that set there are countably many chain-recurrence classes?
\end{quest}

The motivation for posing this question is the well known Smale's
conjecture asserting that for any surface, there is a $C^1$-open
and dense subset of diffeomorphisms of the surface which are
hyperbolic (and in particular, they are tame). Since the dynamics
in the $C^1$-open set $\cU$ given by Theorem A has all of its
chain-recurrence classes except from one contained in periodic
normally hyperbolic surfaces, it seems that a positive answer to
Smale's conjecture would imply that for $C^1$-generic
diffeomorphisms in $\cU$ there are countably many chain-recurrence
classes.

\subsubsection{} We would like to comment that the techniques here are not enough
to treat the examples given by \cite{BLY}. In fact, we mentioned
that C. Bonatti and K. Shinohara have announced that the
quasi-attractor in the example from \cite{BLY} is accumulated by
homoclinic classes which are not contained in periodic surfaces so
the question about the existence of Milnor attractors in those
examples is not settled by now. In more generality one could ask:

 \begin{quest} Does a $C^1$-generic diffeomorphism admit a Milnor attractor?
 \end{quest}

\subsubsection{} Finally, we would like to pose yet another question regarding the example given in Theorem B.

\begin{quest} Can we say anything about the existence and finiteness of SRB measures and/or entropy maximizing measures for the example of Theorem B? \end{quest}

In view of recent results (\cite{VY} and \cite{RHRHUT}) one could expect that these measures may exist but not be unique.

\subsection{Organization of the paper}
In section \ref{SectionMecanismo} we present a general mechanism
for localizing the chain-recurrence classes different from a given
one as lower dimensional classes. We apply this mechanism in
section \ref{SectionTeoremaA} to prove Theorem A. In section
\ref{SectionTeoremaB} we indicate the differences of the proof of
Theorem B and Theorem A. In appendix A we modify the construction
of Theorem B to show how to construct in every manifold a robustly
transitive attractor with partially hyperbolic splitting
$E^{cs}\oplus E^u$ where $E^{cs}$ is not decomposable as two
$Df$-invariant bundles.

\smallskip
{\bf Acknowledgements}\textit{ I would like to thank Sylvain
Crovisier for his patience, corrections and dedication, and, in
particular, for suggesting the use of the semiconjugacy with the
Anosov maps to study this kind of examples. C. Bonatti, J. Buzzi,
L. Diaz, N. Gourmelon, M. Sambarino and K.Shinohara all kindly
listened to the construction and showed interest in it, J.Buzzi
and M. Sambarino were also helpful in the writing of the paper.
Finally, I would like to thank the referees for many important
suggestions to improve considerably the presentation.}

\section{A mechanism for localizing chain-recurrence classes}\label{SectionMecanismo}

Given a homeomorphism $g:\Gamma \to \Gamma$ where $\Gamma$ is a
compact metric space, we say that $g$ is \emph{expansive} if there
exists $\alpha>0$ such that for any pair of distinct points $x
\neq y \in \Gamma$ there exists $n\in \ZZ$ such that
$d(g^n(x),g^n(y))\geq \alpha$.

\begin{prop}\label{ProposicionMecanismo}
Let $f$ be a $C^1$-diffeomorphism and $U$ a filtrating set such
that its maximal invariant set $\Lambda$ admits a partially
hyperbolic structure $T_\Lambda M = E^{cs} \oplus E^u$ such that
$E^{cs}$ is coherent. Assume that there exists a continuous
surjective map $h: \Lambda \to \Gamma$ and a homeomorphism $g :
\Gamma \to \Gamma$ such that:
\begin{itemize}
\item[-] $g\circ h = h \circ f$. \item[-] $h$ is injective in
unstable manifolds. \item[-] There exists a chain recurrence class
$Q$ such that $h^{-1}(\{h(x)\})$ is contained in $\cF^{cs}_U(x)$
and its topological frontier relative to $\cF^{cs}_U(x)$ is
contained $Q$. In particular $h(Q)=\Gamma$. \item[-] The fibers
$h^{-1}(\{y\})$ are invariant under unstable holonomy. \item[-]
$g$ is expansive.
\end{itemize}
Then, every chain-recurrence class in $U$ different from $Q$ is contained in the preimage of a periodic orbit by $h$.
\end{prop}

For simplicity, the reader can follow the proof assuming that $g$
is an Anosov diffeomorphism we shall make some footnotes when some
differences (which are quite small) appear.

\dem{\!\!\!} Let $R \neq Q$ be a chain recurrence class of $f$.
Then, since $\partial h^{-1}(\{y\}) \en Q$ for every $y\in
\Gamma$, we have that $R \cap int(h^{-1}(\{y\})) \neq \emptyset$
for some $y \in \Gamma$.

Conley's theory gives us an open neighborhood $V$ of $R$ whose
closure is disjoint from $Q$ and such that every two points $x,z
\in R$ are joined by arbitrarily small pseudo-orbits contained in
$V$.

\begin{figure}[ht]
\begin{center}
\input{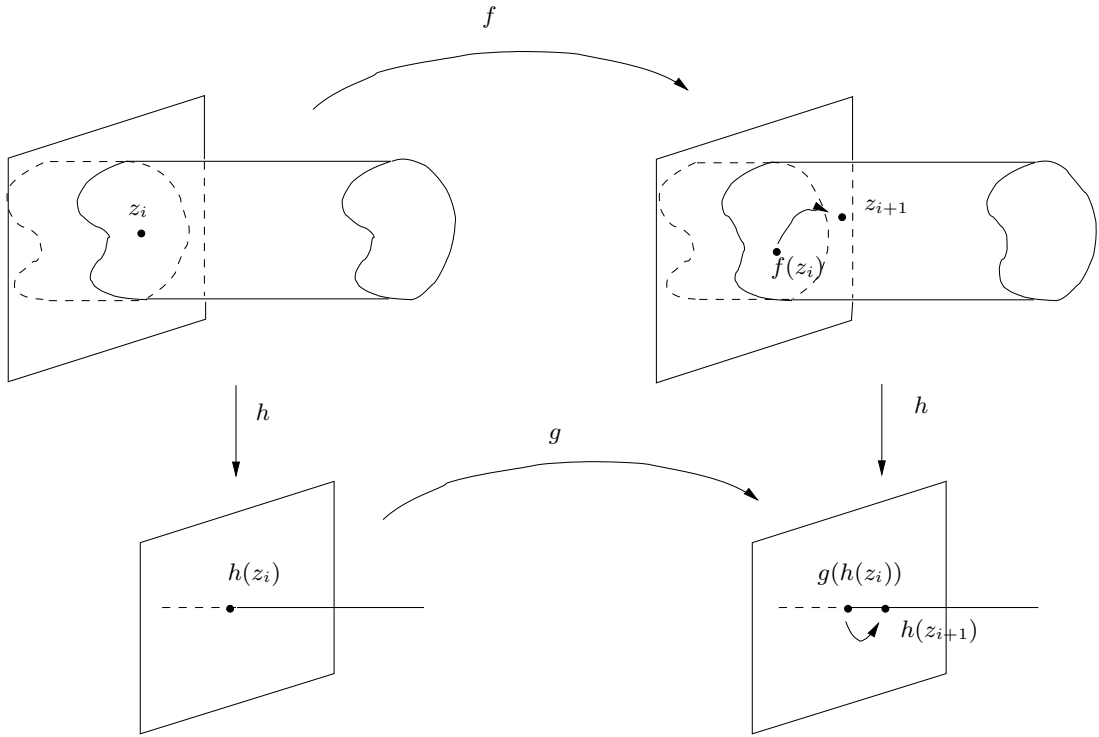}
\caption{\small{Pseudo-orbits for $f$ are sent to pseudo-orbits of $g$ with jumps in the unstable sets.}} \label{FiguraPseudoOrbitas}
\end{center}
\end{figure}

Since $\overline{V}$ does not intersect $Q$, using the invariance
under unstable holonomy of the fibers, we get that there exists
$\eta_0$ such that if $d(x,z)<\eta_0$ and $x\in V$, then $h(x)$
and $h(z)$ lie in the same local unstable manifold\footnote{The
$\zeta-$local unstable set of a point $x$ for an expansive
homeomorphism $g$ is the set of points whose orbit remains at
distance smaller than $\zeta$ for every past iterate. For an
expansive homeomorphism, this set is contained in the unstable
set.}  (it suffices to choose $\eta_0< d(\overline V, Q)$).

Given $\zeta>0$ we choose $\eta>0$ such that $d(x,z)<\eta$ implies
$d(h(x), h(z))< \zeta$. The semiconjugacy implies then that if
$z_0, \ldots z_n$ is a $\eta-$pseudo orbit for $f$, then $h(z_0),
\ldots, h(z_n)$ is a $\zeta$-pseudo orbit for $g$ (that is,
$d(g(h(z_i)), h(z_{i+1})) < \zeta$). Also, if $\eta<\eta_0$ and
$z_0, \ldots z_n$ is contained in $V$, then we get that the the
pseudo-orbit $h(z_0), \ldots, h(z_n)$ has jumps inside local
unstable sets (i.e. $h(z_{i+1}) \in W^{u}_{\zeta}(g(h(z_i)))$).

Take $x \in R$. Then, for every $\eta < \eta_0$  we take $x=z_0,
z_1, \ldots, z_n=x$ ($n\geq 1$) a $\eta-$pseudo orbit contained in
$V$ joining $x$ to itself. Thus, we have that

$$g^n(W^{u}(h(x)))= W^u(h(x))$$

\noindent so, $W^u(h(x))$ is the unstable manifold for $g$ of a
periodic orbit $\mathcal O$. Since $R$ is $f$-invariant and since
the semiconjugacy implies that $f^{-n}(x)$ accumulates on
$h^{-1}(\cO)$, we get that $R$ intersects the fiber
$h^{-1}(\mathcal O)$.

We must now prove that $R \en h^{-1}(\cO)$ which concludes.

Given $\eps>0$ there exists $\delta> 0$ such that if $z_0, \ldots
z_n$ is a $\delta-$pseudo orbit for $g$ with jumps in the unstable
manifold, then $z_n \in W^u_{\eps}(\cO)$ implies that $z_0 \in
W^u_{2\eps}(\cO)$ (notice that a pseudo orbit with jumps in the
unstable manifold of a periodic orbits can be regarded as a pseudo
orbit for a homothety\footnote{In the general case of $g$ being an
expansive homeomorphism, it is very similar since one has that
restricted to the unstable set of a periodic orbit, one can obtain
a metric inducing the same topology where $g^{-1}$ is an uniform
contraction. This follows from \cite{F} and can also be deduced
using the uniform expansion of $f$ in unstable leaves and the
injectivity of the semi-conjugacy along unstable leaves.} in
$\RR^k$).

Assume that there is a point $z\in R$ such that $h(z) \in
W^u(\mathcal O) \backslash \mathcal O$. So, there are arbitrarily
small pseudo orbits contained in $V$ joining $z$ with a point  in
$h^{-1}(\mathcal O)$. This implies that after sending the pseudo
orbit by $h$ we would get arbitrarily small pseudo orbits for $g$,
with jumps in the unstable manifold, joining $h(z)$ with $\mathcal
O$. This contradicts the remark made in the last paragraph.

So, we get that $R$ is contained in $h^{-1}(\mathcal O)$ where
$\mathcal O$ is a periodic orbit of $g$. \lqqd

\section{Examples in $\TT^3$. Proof of Theorem A}\label{SectionTeoremaA}

\subsection{Construction of the example}\label{SectionConstruccion1}

In this section we shall construct an open set $\cU$ of
$\Diff^r(\TT^3)$ for $r\geq 1$ verifying Theorem A.

The construction is very similar to the one of Carvalho's example
(\cite{Car}) following \cite{BV} with the difference that instead
of creating a source, we create an expanding saddle.

\subsubsection{} We start with a linear Anosov diffeomorphism $A:\T^3 \to \T^3$
admitting a splitting $E^s\oplus E^u$ where $\dim E^s=2$.

We assume that $A$ has complex eigenvalues on the $E^s$ direction
so that $E^s$ cannot split as a dominated sum of other two
subspaces. For example, the matrix

\[      \left(
            \begin{array}{ccc}
              1 & 1 & 0 \\
              0 & 0 & 1 \\
              1 & 0 & 0 \\
            \end{array}
          \right)
\]

\noindent which has characteristic polynomial $1 + \lambda^2
-\lambda^3$ works since it has only one real root, and it is
larger than one.

Considering an iterate, we may assume that there exists $\lambda<1/3$ satisfying:

 \[ \|(DA)_{/E^s}\|< \lambda \ \ ; \ \ \|(DA)_{/E^u}^{-1}\|<\lambda  \]

\subsubsection{} Let $q$ and $r$ be different fixed points of $A$.

Consider $\delta$ small enough such that $B(q,6\delta)$ and
$B(r,6\delta)$ are pairwise disjoint and at distance larger than
$400\delta$ (this implies in particular that the diameter of
$\TT^3$ is larger than $400\delta$).


Let  $\cE^u$ be a family of closed cones around the subspace $E^u$
of $A$ which is preserved by $DA$ (that is $D_xA(\cE^u(x))\en int
(\cE^u(Ax))$). We shall consider the cones are narrow enough so
that any curve tangent to $\cE^u$ of length bigger than $L$
intersects any stable disk of radius $\delta$. Let $\cE^{cs}$ be a
family of closed cones around $E^s$ preserved by $DA$.

From now on, $\delta$ remains fixed. Given $\eps>0$ such that
$\eps \ll \delta$,\footnote{If $K$ bounds $\|A\|$ and
$\|A^{-1}\|^{-1}$ then $\frac{\delta}{10K}$ is enough.} we can
choose $\nu$ sufficiently small such that every diffeomorphism $g$
which is $\nu$-$C^0$-close to $A$ is semiconjugated to $A$ with a
continuous surjection $h$ which is $\eps$-$C^0$-close to the
identity (this is a classical result on topological stability of
Anosov diffeomorphisms, see \cite{W}).

\subsubsection{} We shall modify $A$ inside $B(q,\delta)$ such that we get a new
diffeomorphism $F:\T^3 \to \T^3$ that verifies the following
properties:

\begin{itemize}
\item[-] $F$ coincides with $A$ outside $B(q,\delta)$ and lies at $C^0$-distance smaller than $\nu$ from $A$.
\item[-] The point $q$ is a hyperbolic saddle fixed point of stable index $1$ and such that the product of its two eigenvalues with smaller modulus is larger than $1$. We also assume that the length of the stable manifold of $q$ is larger than $\delta$.
\item[-]$D_xF(\cE^u(x)) \en int(\cE^{u}(F(x)))$. Also, for every $w\in \cE^u(x)\setminus \{0\}$ we have $\|DF_x^{-1}w\|< \lambda \|w\|$.
\item[-] $F$ preserves the stable foliation of $A$. Notice that the foliation will no longer be stable.
\item[-] For some small $\beta>0$ we have that $\|D_xF v\| < (1+\beta) \|v\|$ for every $v$ tangent to the stable foliation of $A$ preserved by $F$ and every $x$.
\end{itemize}

\begin{figure}[ht]
\begin{center}
\input{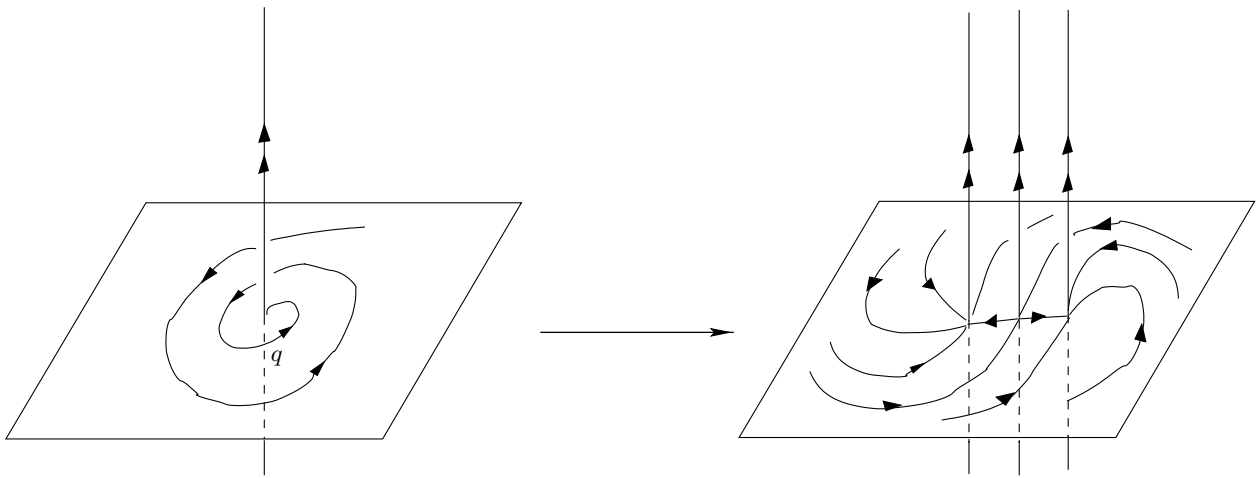}
\caption{\small{Modification of $A$ in a neighborhood of $q$.}} \label{FiguraPseudoOrbitas}
\end{center}
\end{figure}

This construction can be made using classical methods (see
\cite{BV} section 6). Indeed, consider a small neighborhood $U$ of
$q$ such that $U\en B(q,\nu/2)$ such that $U$ admits a chart
$\varphi: U \to \DD^2 \times [-1,1]$ which sends $q$ to $(0,0)$
and sends stable manifolds of $A$ in sets of the form $\DD^2
\times \{t\}$ and unstable ones into sets of the form $\{s\}
\times [-1,1]$. We can modify $A$ by isotopy inside $U$ in such a
way that the sets $\DD^2\times \{t\}$ remain an invariant
foliation but such that the derivative of $q$ becomes the identity
in the tangent space to $\varphi^{-1}(\DD^2\times \{0\})$ which is
invariant and such that the dynamics remains conjugated to the
initial one. At this point, the norm of the images of unit vectors
tangent to the stable foliation of $A$ are not expanded by the
derivative.

Now, one can modify slightly the dynamics in $\varphi^{-1}(\DD^2
\times \{0\})$ in order to obtain the desired conditions on the
eigenvalues of $q$ for $F$. It is not hard to see that for
backward iterates there will be points outside
$\varphi^{-1}(\DD^2\times \{0\})$ which will approach $q$ so one
can obtain the desired length of the stable manifold of $q$ by
maybe performing yet another small modification. All this can be
made in order that the vectors tangent to the stable foliation of
$A$ are expanded by $DF$ by a factor of at most $(1+\beta)$ with
$\beta$ as small as we desire.

The fact that we can keep narrow cones invariant under $DF$ seems
difficult to obtain in view that we made all this modifications.
However, the argument of \cite{BV} (page 190) allows to obtain it:
This is achieved by conjugating the modification with appropriate
homotheties in the stable direction.

The last condition on the norm of $DF$ in the tangent space to the
stable foliation of $A$ seems quite restrictive, more indeed in
view of the condition on the eigenvalues of $q$. This condition
(as well as property (P7) below) shall be only used (and will be
essential) to obtain the ergodic properties of the diffeomorphisms
in the open set we shall construct. Nevertheless, one can
construct such a diffeomorphism as explained above.

\subsubsection{} There exists a $C^1$-open neighborhood $\cU_1$ of $F$ such that for every $f \in \cU_1$ we have that:

\begin{itemize}
\item[(P1)] There exists a continuation $q_f$ of $q$ and $r_f$ of $r$. The point $r_f$ has stable index $2$ and complex eigenvalues. The point $q_f$ is a saddle fixed point of stable index $1$, such that the product of its two eigenvalues with smaller modulus is larger than $1$ and such that the length of the stable manifold is larger than $\delta$.

\item[(P2)] $D_xf(\cE^u(x)) \en int(\cE^{u}(f(x)))$. Also, for
every $w\in \cE^u(x)$ we have $$\|Df_x w\| \geq \lambda^{-1} \|w\|
.$$

\item[(P3)] $f$ preserves a foliation $\cF^{cs}$ which is $C^0$-close to the stable foliation of $A$. Also, each leaf of $\cF^{cs}$ is $C^1$-close to a leaf of the stable foliation of $A$.

\item[(P4)]  For every $x\notin B(q,\delta)$ we have that if $v \in \cE^{cs}(x)$  then $$\|D_x f v\| \leq \lambda \|v\| .$$ This is satisfied for $F$ since $F=A$ outside $B(q,\delta)$.

\item[(P5)] There exists a continuous and surjective map $h_f: \TT^3\to \TT^3$ such that $$h_f \circ f = A\circ h_f$$ \noindent and $d(h(x), x) < \eps$ for every $x\in \TT^3$.
\end{itemize}

The fact that properties (P1), (P2) and (P4) are $C^1-$robust is immediate, robustness of (P5) follows from the choice of $\nu$.

Property (P3) holds in a neighborhood of $F$ since $F$ preserves the stable foliation of $A$ which is a $C^1-$foliation (see \cite{HPS} chapter 7). The foliation $\cF^{cs}$ will be tangent to $E^{cs}$ a bidimensional bundle which is $f$-invariant and contained in $\cE^{cs}$. Other way to proceed in order to obtain an invariant foliation is to use Theorem 3.1 of \cite{BF} of which all hypothesis are verified here but we shall not state it.

Since the cones $\cE^u$ are narrow and from (P3) one has that:

\begin{itemize}
\item[(P6)] Every curve of length $L$ tangent to $\cE^u$ will intersect any disc of radius $2\delta$ in $\cF^{cs}$.
\end{itemize}

Finally, there exists an open set $\cU_2 \en \cU_1$ such that for $f\in \cU_2$ we have:

\begin{itemize}
\item[(P7)]   $\|D_xf v\| \leq (1+\beta) \|v\|$ for every $v\in \cE^{cs}(x)$ and every $x$.
\end{itemize}

\subsubsection{} We shall close this section by proving that for these examples there exists a unique quasi-attractor for the dynamics.

\begin{lema}\label{unicocasiatractor}
For every $f\in \U_1$ there exists an unique quasi-attractor $Q_f$. This quasi attractor contains the homoclinic class of $r_f$, the continuation of $r$.
\end{lema}

\dem{\!\!} We use the same argument as in \cite{BLY}.

 There is a center stable disc of radius bigger than $2\delta$ contained in the stable manifold of $r_f$ ((P3) and (P4)). So, every unstable manifold of length bigger than $L$ will intersect the stable manifold of $r_f$ ((P6)).

Let $Q$ be a quasi attractor, so, there exists a sequence  $U_n$, of neighborhoods of $Q$ such that $f(\overline{ U_n}) \en U_n$ and $Q=\bigcap_n \overline{U_n}$.

Since $U_n$ is open, there is a small unstable curve $\gamma$ contained in $U_n$. Since $Df$ expands vectors in $\cE^u$ we have that the length of $f^k(\gamma)$ tends to $+\infty$ as $n\to +\infty$. So, there exists $k_0$ such that $f^{k_0}(\gamma) \cap W^s(r_f) \neq \emptyset$. So, since $f(\overline{U_n})\en U_n$ we get that $U_n \cap W^s(r_f) \neq \emptyset$, using again the forward invariance of $U_n$ we get that $r_f \in \overline{U_n}$.

This holds for every $n$ so $r_f \in Q$. Since the homoclinic class of $r_f$ is chain transitive, we also get that $H(r_f)\en Q$.

From Conley's theory (cf. \ref{DefinicionCONLEY}), every homeomorphism of a compact metric space there is at least one chain recurrent class which is a quasi attractor. This concludes.
\lqqd

\subsection{The example verifies the mechanism}
We shall consider $f\in \cU_1$ so that it verifies (P1)-(P6).

\subsubsection{} Let $\cA^s$ and $\cA^u$ be, respectively, the stable and unstable foliations of $A$, which are linear foliations. Since $A$ is a linear Anosov diffeomorphism, the distances inside the leaves of the foliations and the distances in the manifold are equal in small neighborhoods of the points if we choose a convenient metric.

Let $\cA^s_\eta(x)$ denote the ball of radius $\eta$ around $x$ inside the leaf of $x$ of $\cA^s$. For any $\eta>0$, it is satisfied that $A(\cA^s_\eta(x))\en \cA^s_{\eta/3}(Ax)$ (an analogous property is satisfied by $\cA^u_\eta(x)$ and backward iterates).

\subsubsection{} The distance inside the leaves of $\Fol^{cs}$ is similar to the ones in the ambient manifold since each leaf of $\cF^{cs}$ is $C^1$-close to a leaf of $\cA^s$. That is, there exists $\rho \approx 1$ such that if $x,y$ belong to a connected component of $\Fol^{cs}(z) \cap B(z,10\delta)$ then $\rho^{-1}d_{cs}(x,y)< d(x,y) <\rho d_{cs}(x,y)$ where $\Fol^{cs}(z)$ denotes the leaf of the foliation passing through $z$ and $d_{cs}$ the distance restricted to the leaf.

For $z\in \TT^3$ we define $W^{cs}_{loc}(z)$ (the \emph{local center stable manifold} of $z$) as the $2\delta$-neighborhood of $z$ in $\cF^{cs}(z)$ with the distance $d_{cs}$.

Also, we can assume that for some $\gamma < \min \{\|A\|^{-1},\|A^{-1}\|^{-1}, \delta/10\}$ the plaque $W^{cs}_{loc}(x)$ is contained in a $\gamma/2$ neighborhood of $\cA^{s}_{2\delta}(x)$, the disc of radius $2\delta$ of the stable foliation of $A$ around $x$.

\begin{lema}\label{trapping} We have that $f(\overline{W^{cs}_{loc}(x)}) \en W^{cs}_{loc}(f(x))$.
\end{lema}

\dem{\!\!} Consider around each $x\in \T^3$ a continuous map $b_x:\D^2 \times [-1,1] \to \T^3$ such that $b_x(\{0\} \times [-1,1])=\cA^u_{3\delta}(x)$ and $b_x(\D^2 \times \{t\})= \cA^s_{3\delta}(b_x(\{0\}\times\{t\}))$. For example, one can choose $b_x$ to be affine in each coordinate to the covering of $\TT^3$.

  Thus, it is not hard to see that one can assume also that $b_x(\frac 1 3 \D^2 \times \{t\} ) = \cA^s_{\delta}(b_x(\{0\}\times\{t\}))$ and that $b_x(\{y\} \times [-1/3,1/3])= \cA^u_{\delta}(b_x(\{y\}\times \{0\}))$. Let

  $$B_x = b_x(\D^2 \times [-\gamma/2,\gamma/2]) .$$

We have that $A(B_x)$ is contained in $b_{Ax}( \frac{1}{3}\D^2 \times [-1/2,1/2])$. Since $f$ is $\eps$-$C^0$-near $A$, we get that $f(B_x) \en b_{f(x)}( \frac{1}{2}\D^2 \times [-1,1])$.

Let $\pi_1: \D^2 \times [-1,1] \to \D^2$ such that $\pi_1(x,t)=x$. We have that $\pi_1(b_{f(x)}^{-1}(W^{cs}_{loc}(f(x))))$ contains $\frac{1}{2}\D^2$ from how we chose $\gamma$ and from how we have defined the local center stable manifolds\footnote{In fact, $b_{f(x)}^{-1}(W^{cs}_{loc}(g(x))) \cap \frac{1}{2}\D^2\times [-1,1]$ is the graph of a $C^1$ function from $\frac{1}{2}\D^2$ to $[-\gamma/2,\gamma/2]$ if $b_x$ is well chosen.}.

Since $f(\Fol^{cs}(x)) \en \Fol^{cs}(f(x))$ and $f(W^{cs}_{loc}(x)) \en  b_{f(x)}( \frac{1}{2}\D^2 \times [-1,1])$ we get the desired property.
\lqqd

\subsubsection{} The fact that $f\in \U_1$ is semiconjugated with $A$ together with the fact that the semiconjugacy is $\eps$-$C^0$-close to the identity gives us the following easy properties about the fibers (preimages under $h_f$) of the points.

As in \ref{DefinicionPartialHyp}, we denote

$$\Pi^{uu}_{x,z}: U\en W^{cs}_{loc}(x) \to W^{cs}_{loc}(z)$$

\noindent the unstable holonomy where $z\in \cF^{u}(x)$ and $U$ is a neighborhood of $x$ in $W^{cs}_{loc}(x)$ which can be considered large if $z$ is close to $x$ in $\cF^u(x)$. In particular, let $\gamma>0$ be such that if $z \in \cF^{u}_{\gamma}(x)$ then the holonomy is defined in a neighborhood of radius $\eps$ of $x$.

\begin{prop}\label{propiedades} Consider $y= h_f(x)$ for $x\in \TT^3$:
\begin{enumerate}
\item $h_f^{-1}(\{y\})$ is a compact connected set contained in $W^{cs}_{loc}(x)$.
\item If $z\in \cF^{u}_{\gamma}(x)$, then $h_f(\Pi^{uu}_{x,z}(h^{-1}_f (\{y\})))$ is exactly one point.
\end{enumerate}
\end{prop}
\dem{\!\!} (1)  Since $h_f$ is $\eps$-$C^0$-close the identity, we get that for every point $y\in \T^3$, $h_f^{-1}(\{y\})$ has diameter smaller than $\eps$. Since $\eps$ is small compared to $\delta$, it is enough to prove that $h_f^{-1}(\{y\}) \en W^{cs}_{loc}(x)$ for some $x\in h_f^{-1}(\{y\})$.

Assume that for some $y\in \T^3$, $h_f^{-1}(\{y\})$ intersects two different center stable leaves of $\Fol^{cs}$ in points $x_1$ and $x_2$.

Since the points are near, we have that $\cF^{u}_{\gamma}(x_1) \cap W^{cs}_{loc}(x_2) =\{z\}$. Thus, by forward iteration, we get that for some $n_0>0$ we have $d(f^{n_0}(x_1), f^{n_0}(z))> 3\delta$.

Lemma \ref{trapping} gives us that $d(f^{n_0}(x_2), f^{n_0}(z)) < 2\delta$ and so, we get that $d(f^{n_0}(x_1),f^{n_0}(x_2)) > \delta$ which is a contradiction since $\{f^{n_0}(x_1),f^{n_0}(x_2)\} \en h_f^{-1}(\{A^{n_0}(y)\})$ which has diameter smaller than $\eps \ll \delta$.

Also, since the dynamics is trapped in center stable manifolds, we get that the fibers must be connected since one can write them as

$$ h^{-1}(\{h(x)\})= \bigcap_{n\geq 0} f^n(W^{cs}_{loc}(f^{-n}(x))) .$$

(2)  Since $f^{-n}(h_f^{-1}(\{y\})) = h_f^{-1}(\{A^{-n}(y)\})$ we get that $\diametro(f^{-n}(h_f^{-1}(\{y\}))) < \eps$ for every $n>0$.

This implies that there exists $n_0$ such that if $n>n_0$ then $f^{-n}(\Pi^{uu}_{x,z}(h^{-1}_f (\{y\})))$ is sufficiently near $f^{-n}(h_f^{-1}(\{y\}))$. So, we have that $$\diametro(f^{-n}(\Pi_{x,z}^{uu}(h_f^{-1}(\{y\})))) < 2\eps \ll \delta .$$

Assume that $h_f(\Pi_{x,z}^{uu}(h_f^{-1}(\{y\})))$ contains more than one point. These points must differ in the stable coordinate of $A$, so, after backwards iteration we get that they are at distance bigger than $3\delta$. Since $h_f$ is $\eps$-$C^0$-close the identity this represents a contradiction.
 \lqqd

\begin{obs} The second statement of the previous proposition gives that the fibers of $h_f$ are invariant under unstable holonomy.\finobs
\end{obs}

\subsubsection{} The following simple lemma is essential in order to satisfy the properties of Proposition \ref{ProposicionMecanismo}.

\begin{lema}\label{conexosgrandes} For every $f\in \U_1$, given a disc $D$ in $W^{cs}_{loc}(x)$ whose image by $h_f$ has at least two points, then $D \cap \cF^{u}(r_f) \neq \emptyset$ and the intersection is transversal.
\end{lema}

\dem{\!\!} Given a subset $K\en \Fol^{cs}(x)$ we define its \emph{center stable diameter} as the diameter with the metric $d_{cs}$ defined above induced by the metric in the manifold. We shall first prove that there exists $n_0$ such that $\diametro_{cs} (f^{-n_0}(D)) > 100 \delta$:

Since $D$ is arc connected so is $h_f(D)$, so, it is enough to suppose that $\diametro (D)<\delta$.
We shall first prove that $h_f(D)$ is contained in a stable leaf of the stable foliation of $A$. Otherwise, there would exist points in $h_f(D)$ whose future iterates separate more than $2\delta$, this contradicts that the center stable plaques are trapped for $f$ (Lemma \ref{trapping}).

One now has that, since $A$ is Anosov and that $h_f(D)$ is a connected compact set with more than two points contained in a stable leaf of the stable foliation, there exists $n_0 > 0$ such that $A^{-n_0}(h_f(D))$ has stable diameter bigger than $200\delta$ (recall that $\diametro \TT^3 > 400\delta$). Now, since $h_f$ is close to the identity, one gets the desired property.

We conclude by proving the following:

\noindent{\bf Claim 1:} If there exists $n_0$ such that $f^{-n_0}(D)$ has diameter larger than $100\delta$, then $D$ intersects $\cF^{u}(r_f)$.

\medskip
\dem{of the claim} This is proved in detail in section 6.1 of \cite{BV} so we shall only sketch it.

If $f^{-n_0}(D)$ has diameter larger than $100\delta$, from how we choose $\delta$ we have that there is a compact connected subset of $f^{-n_0}(D)$ of diameter larger than $35\delta$ which is outside $B(q,6\delta)$.

So, $f^{-n_0-1}(D)$ will have diameter larger than $100\delta$ and the same will happen again. This allows to find a point $x\in D$ such that $\forall n>n_0$ we have that $f^{-n}(x)\notin B(q,6\delta)$.

Now, considering a small disc around $x$ we have that by backward iterates it will contain discs of radius each time bigger and this will continue while the disc does not intersect $B(q,\delta)$. If that happens, since $f^{-n}(x) \notin B(q,6\delta)$ the disc must have radius at least $3\delta$.

This proves that there exists $m$ such that $f^{-m}(D)$ contains a center stable disc of radius bigger than $2\delta$, so, the unstable manifold of $r_f$ intersects it. Since the unstable manifold of $r_f$ is invariant, we deduce that it intersects $D$ and this concludes the proof of the claim.

Transversality of the intersection is immediate from the fact that $D$ is contained in $\cF^{cs}$ which is transversal to $\cF^u$.
\lqqd

\subsubsection{} We obtain the following corollary which puts us in the hypothesis of Proposition \ref{ProposicionMecanismo}:

\begin{cor}\label{interior} For every $f\in \U_1$, let $x\in \partial h_f^{-1}(\{y\})$ $($relative to the local center stable manifold of  $h_f^{-1}(\{y\}))$, then, $x$ belongs to the homoclinic class of $r_f$, and in particular, to $Q_f$.
\end{cor}

\dem{\!\!} Notice first that the stable manifold of $r_f$ coincides with $\cF^{cs}(r_f)$ which is dense in $\TT^3$. This follows from the fact that when iterating an unstable curve, it will eventually intersect the stable manifold of $r_f$, since the stable manifold of $r_f$ is invariant, we obtain the density of $\cF^{cs}(r_f)$.

Now, considering $x\in \partial h_f^{-1}(\{y\})$, and $\eps>0$, we consider a connected component $\tilde D$ of $\cF^{cs}(r_f) \cap B(x,\eps)$. Clearly, since the fibers are invariant under holonomy and $x \in \partial h_f^{-1}(\{y\})$ we get that $\tilde D$ contains a disk $D$ which is sent by $h_f$ to a non trivial connected set. Using the previous lemma we obtain that there is a homoclinic point of $r_f$ inside $B(x,\eps)$ which concludes.
\lqqd

\subsubsection{} The following corollary will allow us to use Theorems \ref{BDP} and \ref{BDPV}.

\begin{cor}\label{qgenLambdag} For every $f\in \U_1$ we have that $q_f\in H(r_f)$.
\end{cor}

\dem{\!\!} Consider $U$, a neighborhood of $q_f$, and $D$ a center stable disc contained in $U$.

Since the stable manifold of $q_f$ has length bigger than $\delta>\eps$, after backward iteration of $D$ one gets that $f^{-k}(D)$ will eventually have diameter larger than $\eps$, thus $h_f(D)$ will have at least two points, this means that $q_f \in \partial h_f^{-1}(\{h(q_f)\})$. Corollary \ref{interior} concludes.
\lqqd

\subsubsection{}\label{SectionTeorema37} We finish this section by proving the following theorem which is the topological part of Theorem A.

\begin{teo}\label{TeoParteTopologica}
\begin{itemize}
    \item[(i)] For every $f\in \cU_1$ there exists a unique quasi-attractor $Q_f$ which contains the homoclinic class $H(r_g)$ and such that every chain-recurrence class $R\neq Q_f$ is contained in a periodic disc of $\cF^{cs}$.
    \item[(ii)] For every $f \in \cG_{BC} \cap \cG_{BDP} \cap \cU_1$ we have that $H(r_f)=Q_f$ and is contained in the closure of the sources of $f$.
    \item[(iii)] For every $r\geq 2$, there exists a $C^1$-open dense subset $\cU_3$ of $\cU_1$ and a $C^r$-residual subset $\cG^r \en \cU_3 \cap \Diff^r(\TT^3)$ such that for every $f\in \cG^r$ the homoclinic class $H(r_f)$ intersects the closure of the sources of $f$.
    \item[(iv)] For every $f\in \cU_1$ there exists a unique Milnor attractor contained in $Q_f$.
\end{itemize}
\end{teo}

\dem{\!\!} Part (i) follows from Proposition \ref{ProposicionMecanismo} since $h_f$ is the desired semiconjugacy: Indeed, Proposition \ref{propiedades} and Corollary \ref{interior} show that the hypothesis of the mentioned proposition are verified (notice that $A$ is clearly expansive).

Part (ii) follows from Theorem \ref{BDP} using Corollary \ref{qgenLambdag}. Notice that $E^{cs}$ cannot be decomposed in two $Df$-invariant subbundles since $Df$ has complex eigenvalues in $r_f$.

Similarly, part (iii) follows from Theorem \ref{BDPV} and Corollary \ref{BDPV}. The need for considering $\cU_3$ comes from \cite{BD} (see Theorem \ref{BDPV}).

To prove (iv) notice that every point which does not belong to the fiber of a periodic orbit belongs to the basin of $Q_f$: Since there are only countably many periodic orbits and their fibers are contained in two dimensional discs (which have zero Lebesgue measure) this implies directly that the basin of $Q_f$ has total Lebesgue measure:

Consider a point $x$ whose omega-limit set $\omega(x)$ is contained in a chain recurrence class $R$ different from $Q_f$. Then, since this chain recurrence class is contained in the fiber $h_f^{-1}(\cO)$ of a periodic orbit $\cO$ of $A$, which in turn is contained in the local center stable manifold of some point $z\in \TT^3$. This implies that some forward iterate of $x$ is contained in $W^{cs}_{loc}(z)$. The fact that the dynamics in $W^{cs}_{loc}$ is trapping (see Lemma \ref{trapping}) and the fact that $\partial h_f^{-1}(\cO) \en Q_f$ (see Corollary \ref{interior}) gives that $x$ itself is contained in $h_f^{-1}(\cO)$ as claimed.

Now, Lemma 1 of \cite{Milnor} implies that $Q_f$ contains an attractor in the sense of Milnor.
\lqqd

We have just proved parts (a), (b), (c) and (d) of Theorem A hold in $\cU_3$. In fact, for the $C^1$-topology, we have obtain a slightly stronger property than (d) holds in $\cU_1$. Also, we have proved that (e) is satisfied.

\begin{obs} The choice of having complex eigenvalues for $A$ was only used to guaranty that $E^{cs}$ admits no $Df$-invariant subbundles. One could have started with any linear Anosov map $A$ and modify the derivative of a given fixed or periodic point $r$ to have complex eigenvalues and the construction would be the same. \end{obs}

\subsection{Ergodic properties} In this section we shall work with $f\in \cU_2$ so that properties (P1)-(P7) are verified.

\subsubsection{} We shall briefly explain how it can be deduced from \cite{BV} that there exists a unique SRB measure for every $f\in \cU_2$ of class $C^2$. Let us first consider $U$ a small neighborhood of $Q_f$ and

$$\Lambda_f= \bigcap_{n\geq 0} f^n(\overline U)$$

\noindent which is a (not-necessarily transitive) topological attractor.

We shall show that the hypothesis of Theorem A of \cite{BV} are satisfied for $\Lambda_f$ (see also Theorem 11.25 in \cite{BDV}), and thus, we get that there are at most finitely many SRB measures such that the union of their (statistical) basins has full Lebesgue measure in the topological basin of $\Lambda_f$. Clearly, for $\Lambda_f$ one has (H1) and (H2). Hypothesis (H3) follows from the following:

\begin{prop}\label{exponentes} For every $x\in \TT^3$ and $D\en W^{uu}_{loc}(x)$ an unstable arc, we have full measure set of points which have negative Lyapunov exponents in the direction $E^{cs}$.
\end{prop}

\dem{\!\!} The proof is exactly the same as the one in Proposition 6.5 of \cite{BV} so we omit it. Notice that conditions (P2), (P4) and (P7) in our construction imply conditions $(i)$ and $(ii)$ in section 6.3 of \cite{BV}.
\lqqd

\subsubsection{} The set $\Lambda_f$ does not verify the hypothesis of Theorem B of \cite{BV} since we do not have minimality of the unstable foliation.

However, the fact that the stable manifold of $r_f$ contains $W^{cs}_{loc}(r_f)$, gives that every unstable manifold intersects $W^s(r_f)$ and so we get that every compact subset of $\Lambda_f$ saturated by unstable sets must contain $\overline{\cF^u(r_f)}$. This implies that for every $x\in \overline{\cF^u(r_f)}$ we have that $\overline{\cF^u(r_f)} = \overline{\cF^u(x)}$ and $\overline{\cF^u(r_f)}$ is the only compact set with this property (we say that $\overline{\cF^u(r_f)}$ is the unique minimal set of the foliation $\cF^u$).

It is not hard to see how the proof of \cite{BV} works in this context\footnote{See the first paragraph of section 5 in \cite{BV}. Our Proposition \ref{exponentes} implies that (H3) is verified. Moreover, every unstable arc converges after future iteration to the whole $\overline{\cF^{u}(r_g)}$, and since the unstable foliation is minimal in $\overline{\cF^{u}(r_g)}$ we get that there is only one accessibility class there as needed for their Theorem B.}. We get thus that $f$ admits an unique SRB measure $\mu$ and clearly, the support of this SRB measure is $\overline{\cF^{u}(r_f)}$.

We claim that $\overline{\cF^{u}(r_f)}=H(r_f)$: this follows from the fact that the SRB measure $\mu$ is hyperbolic (by Proposition \ref{exponentes}) and that the partially hyperbolic splitting separates the positive and negative exponents of $\mu$ (this is given in Proposition 1.4 of \cite{C} which states that when one has a hyperbolic measure $\mu$ whose supports admits a dominated splitting respecting the exponents of $\mu$ then the support is contained in a homoclinic class).

\subsubsection{} Finally, since the SRB measure has total support and almost every point converges to the whole support, we get that the attractor is in fact a minimal attractor in the sense of Milnor. We have proved:

\begin{prop}\label{SRB} If $f\in \cU_2$ is of class $C^2$, then $f$ admits a unique SRB measure whose support coincides with $\overline{\cF^{u}(r_f)}=H(r_f)$. In particular, $\overline{\cF^{u}(r_f)}$ is a minimal attractor in the sense of Milnor for $f$.
\end{prop}

The importance of considering $f$ of class $C^2$ comes from the fact that with lower regularity, even if we knew that almost every point in the unstable manifold of $r_f$ has stable manifolds, we cannot assure that these cover a positive measure set due to the lack of absolute continuity in the center stable foliation.

\subsubsection{} However, the information we gathered for smooth systems in $\cU_2$ allows us to extend the result for $C^1$-generic diffeomorphisms in $\cU_2$. Recall that for a $C^1$-generic diffeomorphisms $f\in \cU_2$, the homoclinic class of $r_f$ coincides with $Q_f$.

\begin{teo}\label{AtractorMinimalGenerico} There exists a $C^1$-residual subset $\cG_M \en \cU_2$ such that for every $f\in \cG_M$ the set $Q_f =H(r_f)$ is a minimal Milnor attractor.
\end{teo}

\dem{\!\!} Notice that since $r_f$ has a well defined continuation in $\cU_2$, it makes sense to consider the map $f \mapsto \overline{\cF^{u}(r_f)}$ which is naturally semicontinuous with respect to the Haussdorff topology. Thus, it is continuous in a residual subset $\cG_1$ of $\cU_2$. Notice that since the semicontinuity is also valid in the $C^2$-topology, we have that $\cG_1 \cap \Diff^2(\TT^3)$ is also residual in $\cU_2 \cap \Diff^2(\TT^3)$.

It suffices to show that the set of diffeomorphisms in $\cG_1$ for which $\overline{\cF^{u}(r_f)}$ is a minimal Milnor attractor is a $G_\delta$ set (countable intersection of open sets) since we have already shown that $C^2$ diffeomorphisms (which are dense in $\cG_1$) verify this property.

Given an open set $U$, we define

$$U^+(f)= \bigcap_{n\leq 0} f^n(\overline{U}).$$

Let us define the set $\cO_U(\eps)$ as the set of $f\in \cG_1$ such that they satisfy one of the following (disjoint) conditions

\begin{itemize}
\item[-] $\overline{\cF^{u}(r_f)}$ is contained in $U$ \emph{or}
\item[-] $\overline{\cF^{u}(r_f)} \cap \overline{U}^c \neq \emptyset$ and $Leb(U^+(f))<\eps$
\end{itemize}

We must show that these sets are open in $\cG_1$ (it is not hard
to show that if we consider a countable basis of the topology and
$\{U_n\}$ are finite unions of open sets in the basis then
$\cG_M=\bigcap_{n,m} \cO_{U_n}(1/m)$).

To prove that these sets are open, we only have to prove the semicontinuity of the measure of $U^+(f)$ (since the other conditions are clearly open from how we chose $\cG_1$).

Let us consider the set $\tilde K = \overline U \backslash U^+(f)$, so, we can write $\tilde K$ as an increasing union $\tilde K = \bigcup_{n\geq 1} K_n$ where $K_n$ is the set of points which leave $\overline{U}$ in less than $n$ iterates.

So, if $Leb(U^+(f))<\eps$, we can choose $n_0$ such that $Leb(\overline{U} \backslash K_{n_0}) < \eps$, and in fact we can consider $K_{n_0}'$ a compact subset of $K_{n_0}$ such that $Leb(\overline{U}\backslash K'_{n_0})<\eps$.

In a small neighborhood $\cN$ of $f$, we have that if $f'\in \cN$, then $K'_{n_0} \en \overline{U} \backslash U^+(f')$. This concludes.
\lqqd

This completes the proof of part (f) of Theorem A.

\subsubsection{}
To conclude the proof of Theorem A in $\cU = \cU_2 \cap \cU_3$ one
only has to check property (g) which follows directly from the
main theorem of \cite{BF}.

\section{Examples in any manifold. Proof of Theorem B}\label{SectionTeoremaB}

In this section we shall show how to construct an example verifying Theorem B. We shall see that we can construct a quasi-attractor with a partially hyperbolic splitting $E^{cs} \oplus E^u$ such that $E^{cs}$ admits no sub-dominated splitting. In case $E^{cs}$ is volume contracting, it will turn out that this quasi-attractor is in fact a robustly transitive attractor (thus providing examples of robustly transitive attractors with splitting $E^{cs} \oplus E^u$ in every $3-$dimensional manifold) and when there is a periodic saddle of stable index $1$ and such that the product of any two eigenvalues is greater than one and using Theorems \ref{BDP} and \ref{BDPV} we shall obtain that the quasi-attractor will not be isolated for generic diffeomorphisms in a neighborhood.

We shall work only in dimension $3$. It will be clear that by multiplying the examples here with a strong contraction, one can obtain examples in any manifold of any dimension.

A main difference between this construction and the one done in section \ref{SectionTeoremaA} is the use of \emph{blenders} instead of the argument \`a la Bonatti-Viana. Blenders were introduced in \cite{BDPer} and constitute a very powerful tool in order to get robust intersections between stable and unstable manifolds of compact sets. We shall only state some of their properties and not enter in their definition or construction for which there are many excellent references (we recommend chapter 6 of \cite{BDV} in particular).

\subsection{Construction of the example}

\subsubsection{} Let us consider $P: \DD^2 \hookrightarrow \DD^2$ the map given by the Plykin attractor in the disk $\DD^2$ (see \cite{Rob}).

We have that $P(\DD^2) \en int(\DD^2)$, there exist a hyperbolic attractor $\Upsilon \en \DD^2$ and three fixed sources (we can assume this by considering an iterate).

There is a neighborhood $N$ of $\Upsilon$ which is homeomorphic to
the disc with $3$ holes that satisfies that $P(\overline N) \en N$
and

$$\Upsilon= \bigcap_{n\geq 0} P^n(N) . $$

It is well known that given $\eps>0$, one can choose a finite
number of periodic points $s_1, \ldots, s_N$ and $L>0$ such that
if $A= \bigcup_{i=1}^N W^u_L(s_i)$, then, for every $x\in \Upsilon
\setminus A$ one has that $A$ intersects both connected components
of $W^s_\eps (x) \setminus \{x\}$.


We now consider the map $F_0: \DD^2 \times S^1 \hookrightarrow
\DD^2 \times S^1$ given by $F_0(x,t)= (P(x), t)$ whose chain
recurrence classes consist of the set $\Upsilon \times S^1$ which
is a (non transitive) partially hyperbolic attractor and three
repelling circles.

\subsubsection{} In \cite{BDPer} they make a small $C^\infty$ perturbation $F_1$ of
$F_0$, for whom the maximal invariant set in $U= N\times S^1$
becomes a $C^1$-robustly transitive partially hyperbolic attractor
$Q$ which remains homeomorphic to $\Upsilon \times S^1$.

This attractor has a partially hyperbolic structure of the type
$E^s \oplus E^c \oplus E^u$. One can make this example in order
that it fixes the boundary of $\D^2 \times S^1$, this allows to
embed this example (and all the modifications we shall make) in
any isotopy class of diffeomorphisms of any $3$-dimensional
manifold (since every diffeomorphism is isotopic to one which
fixes a ball, then one can introduce this map by a simple
surgery).

\subsubsection{}\label{SectionBlenders} We shall now present $cu$-blenders by their properties: A
$cu$-\emph{blender} $K$ for a diffeomorphism $f: M \to M$ is a
compact $f$-invariant hyperbolic set with splitting $T_K M =
E^{ss} \oplus E^s \oplus E^{u}$ such that the following properties
are verified:

\begin{itemize}
\item[-] $K$ is the maximal invariant subset in a neighborhood
$U$. \item[-] There exists a cone-field $\cE^{ss}$ around $E^{ss}$
defined in all $U$ which is invariant under $Df^{-1}$. \item[-]
There exists a compact region $B$ with non-empty interior (which
is called \emph{activating region}) such that every curve
contained in $U$, tangent to $\cE^{ss}$ with length larger than
$\delta$ and intersecting $B$ verifies that it intersects the
unstable manifold of a point of $K$. \item[-] There exists an open
neighborhood $\cU$ of $f$ such that for every $g$ in $\cU$ the
properties above are verified for the same cone field, the same
set $B$ and for $K_g$ the maximal invariant set of $U$.
\end{itemize}

For more properties and construction of $cu$-blenders, see
\cite{BDV} chapters 6 and \cite{BDPer} (they treat mainly
$cs$-blenders which are $cu$-blenders for $f^{-1}$). There one can
see a proof of the following:

\begin{prop}[\cite{BDPer} Lemma 1.9, \cite{BDV} Lemma 6.8]\label{PropBlender}
If the stable manifold of a periodic point $p\in M$ of stable
index $1$ contains an arc $\gamma$ tangent to $\cE^{ss}$ and
intersecting the activating region of a $cu$-blender $K$, then,
$W^u(p) \en \overline{W^u(q)}$ for every $q$ periodic point in
$K$.
\end{prop}

\subsubsection{} In \cite{BDPer} the diffeomorphism $F_1$ constructed verifies the following properties (see \cite{BDPer} section 4.a page 391, also one can find the indications in \cite{BDV} section 7.1.3):

\begin{itemize}
\item[(F1)] $F_1$ leaves invariant a $C^1$-lamination $\cF^{cs}$
(see \cite{HPS} chapter 7 for a precise definition) tangent to
$E^s \oplus E^c$ whose leaves are homeomorphic to $\RR\times
\SS^1$. \item[(F2)] There are periodic points $p_1, \ldots, p_N$
of stable index $1$ which are homoclinically related and such that
for every $x\in Q$ one has that the connected component of
$\cF^{cs}(x) \setminus \overline{(W^u_L(p_1) \cup \ldots \cup
W^u_L(p_N))}$ containing $x$ has finite volume for every $x\in Q
\setminus \bigcup_{i=1}^N W^u(p_i)_L$. Here $W^u_L(p_i)$ denotes
the neighborhood of $p_i$ in its unstable manifold with the metric
induced by the ambient. \item[(F3)] There is a periodic point $q$
of stable index $2$ contained in a $cu$-blender $K$ such that for
every $0\leq i\leq N$, the stable manifold of $p_i$ intersects the
activatin region of $K$. By Proposition \ref{PropBlender}, the
unstable manifold of $q$ is dense in the union of the unstable
manifolds of $p_i$.

\item[(F4)] The local stable manifold of $q$ intersects every unstable curve of length larger than $L$.
\end{itemize}

Before we continue, we shall make some remarks on the properties. The hypothesis (F1) on the differentiability of the lamination $\cF^{cs}$ will be used in order to apply the results on normal hyperbolicity of \cite{HPS} (chapter 7, Theorem 7.4). It can be seen in \cite{BDPer} that the construction of $F_1$ is made by changing the dynamics in finitely many periodic circles and this can be done without altering the lamination $\cF^{cs}$ which is $C^1$ before modification. This is in fact not necessary; it is possible to apply the barehanded arguments of the proof of Theorem 3.1 of \cite{BF} in order to obtain that for the modifications we shall make, there will exist a lamination tangent to the bundle $E^{cs}$.

Hypothesis (F2) is justified by the fact that the Plykin attractor verifies the same property and the construction of $F_1$ in \cite{BDPer} is made by changing the dynamics in the periodic points by Morse-Smale diffeomorphisms which give rise property (F2) (see section 4.a. of \cite{BDPer}). Notice that by continuous variation of stable and unstable sets, this condition is $C^1$-robust.

Property (F3) is the essence in the construction of \cite{BDPer},
blenders are the main tool for proving the robust transitivity of
these examples. As explained in \ref{SectionBlenders} this is a
$C^1$-open property.

Property (F4) is given by the fact that the local stable manifold
of $q$ can be assumed to be $W^s_{loc}(s) \times \SS^1$ with a
curve removed, where $s\in \Upsilon$ is a periodic point. This is
also a $C^1$-open property.

\begin{figure}[ht]
\begin{center}
\scalebox{.8}{\input{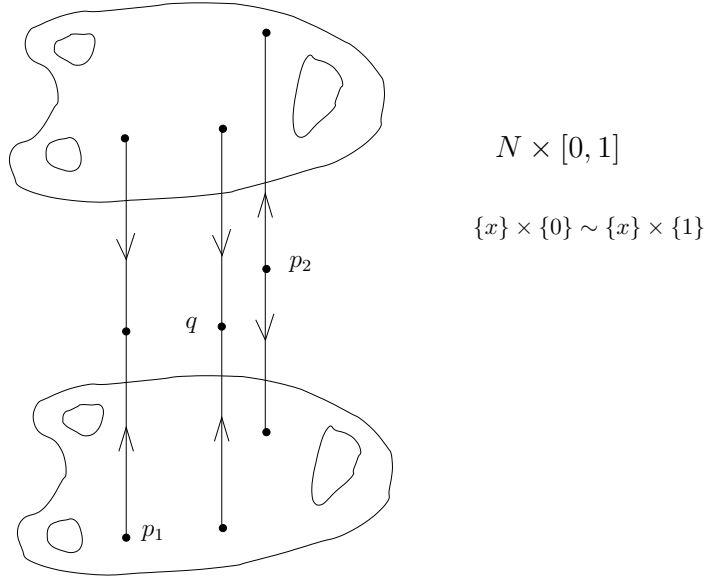}}
\caption{\small{How to construct $F_1$ by small $C^\infty$ perturbations in finitely many circles.}} \label{FiguraAxRodeado}
\end{center}
\end{figure}

\subsubsection{} Let us consider a periodic point $r_1 \in Q$ of stable index $1$
and another one $r_2$ of stable index $2$. We can assume they are
fixed (modulo considering an iterate of $F_1$). Consider
$\delta>0$ small enough such that $B(r_1, 6\delta) \cup B(r_2,
6\delta)$ is disjoint from:

\begin{itemize}
\item[-] the periodic points $p_1, \ldots, p_N, q$ defined above,
\item[-] the blender $K$, \item[-] $\overline{(W^u_L(p_1) \cup
\ldots \cup W^u_L(p_N))}$ and \item[-] from compact connected
pieces of $W^s(p_i)$ ($1\leq i\leq N$) intersecting the activating
region of the blender $K$ and containing $p_i$.
\end{itemize}

\subsubsection{} In the same vein as in section \ref{SectionConstruccion1} we shall construct a diffeomorphism $F_2$ modifying $F_1$ such that:

\begin{itemize}
\item[-] $F_2$ coincides with $F_1$ outside $B(r_1,\delta) \cup
B(r_2,\delta)$. \item[-] $F_2$ preserves the center-stable
lamination of $F_1$. \item[-] $DF_2$ preserves narrow cones
$\cE^u$ and $\cE^{cs}$ around the unstable direction $E^u$ and the
center stable direction $E^{s} \oplus E^c$ of $F_1$ respectively.
Also, vectors in $\cE^{u}$ are expanded uniformly by $DF_2$.
\item[-]The point $r_1$ remains fixed but it becomes a saddle with
stable index $1$ and the product of any pair of eigenvalues of
$r_1$ is larger than $1$. Moreover, the stable manifold of $r_1$
intersects the complement of $B(r_1,6\delta)$. \item[-] The point
$r_2$ remains fixed for $F_2$ but now has complex eigenvalues in
$r_2$.
\end{itemize}

We obtain a $C^1$ neighborhood $\cU_1$ of $F_2$ where for $f\in
\cU$, if we denote

$$\Lambda_f = \bigcap_{n\geq0} f^n(U): $$

\noindent we have:

\begin{itemize}
\item[(P1')] There exists a continuation of the points $p_1,
\ldots, p_N, q, r_1, r_2$ which we shall denote as $p_i(f)$,
$q(f)$ and $r_i(f)$. The point $r_1(f)$ is a saddle of stable
index $1$ and its stable manifold intersects the complement of
$B(r_1,6\delta)$. We further assume that all these periodic points
remain hyperbolic in $\cU$ and that $r_2(f)$ has complex
eigenvalues. \item[(P2')] There is a $Df$-invariant families of
cones $\cE^u$ in $Q_f$ and for every $v\in \cE^u(x)$ we have that
$$\|D_xf v\|\geq \lambda^{-1} \|v\|.$$
\item[(P3')] $f$ preserves a lamination $\cF^{cs}$ which is $C^0$
close to the one preserved by $F_1$ and which is trapped in the
sense that there exists a family $W^{cs}_{loc}(x) \en \cF^{cs}(x)$
such that for every point $x\in Q_f$ the plaque $W^{cs}_{loc}(x)$
is homeomorphic to $(0,1)\times \SS^1$ and verifies that

    $$f(\overline{W^{cs}_{loc}(x)}) \en W^{cs}_{loc}.$$
Moreover, the stable manifold of $r_1(f)$ intersects the
complement of $W^{cs}_{loc}(r_1(f))$.

\item[(P4')] Properties (F2),(F3) and (F4) are satisfied for $f$
and every curve $\gamma$ tangent to $\cE^u$ of length larger than
$L$ intersects the stable manifold of $q(f)$.
\end{itemize}

Notice that (P4') implies that there exists a unique
quasi-attractor $Q_f$ in $U$ for every $f\in \cU$ which contains
the homoclinic class $H(q(f))$ of $q(f)$ (the proof is the same as
Lemma \ref{unicocasiatractor}).

\subsection{The example verifies the mechanism of Proposition \ref{ProposicionMecanismo}}

We shall show that every $f\in \cU$ is in the hypothesis of
Proposition \ref{ProposicionMecanismo} which will conclude the
proof of Theorem B as in Theorem \ref{TeoParteTopologica}.

We shall only sketch the proof since it has the same ingredients
as the proof of Theorem A, the main differences are that instead
of having an a priori semiconjugacy we must construct one and that
we use blenders in order to obtain that the boundary of the fibers
is contained in the unique quasi-attractor.

To construct the semiconjugacy, one uses property (P3'),
specifically the fact that $f(\overline{W^{cs}_{loc}(x)}) \en
W^{cs}_{loc}(x)$ (compare with Lemma \ref{trapping}) to consider
for each point $x\in \Lambda_f$ the set:

$$ A_x = \bigcap_{n\geq 0} f^n(\overline{W^{cs}_{loc}(f^{-n}(x))}) $$

\noindent (compare with Proposition \ref{propiedades} (1)). One
easily checks that the sets $A_x$ constitute a partition of
$\Lambda_f$ into compact connected sets contained in local center
stable manifolds and that the partition is upper-semicontinuous.
It is not hard to prove that if $h_f: \Lambda_f \to
\Lambda_f/_{\sim}$ is the quotient map, then, the map $g:
\Lambda_f /_{\sim} \to \Lambda_f/_{\sim}$ defined such that

$$ h_f \circ f = g \circ h_f $$

\noindent is expansive (in fact, $\Lambda_f/_{\sim}$ can be seen
to be homeomorphic to $\Upsilon$ and $g$ conjugated to $P$). See
\cite{Daverman} for more details on this kind of decompositions
and quotients.

\begin{figure}[ht]
\begin{center}
\input{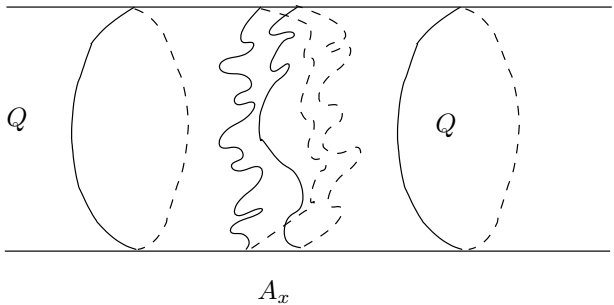}
\caption{\small{The set $A_x$ is surrounded by points in $\overline{W^u(q)} \en Q$ }} \label{FiguraAxRodeado}
\end{center}
\end{figure}

Since fibers are contained in center stable sets, we get that
$h_f$ is injective on unstable manifolds and one can check that
the fibers are invariant under unstable holonomy (see the proof of
Proposition \ref{propiedades} (2)). Stable sets of $g$ are dense
in $\Lambda_f/_{\sim}$.

The point $r_1(f)$ will be in the boundary of
$h^{-1}_f(\{h_f(r(f))\})$ since its stable manifold is not
contained in $W^{cs}_{loc}(r_1(f))$.

\subsubsection{} 
We claim that the boundary of the fibers restricted to center-stable manifolds is contained in the unique quasi-attractor $Q_f$. This is proven as follows:

First, notice that by the $cu$-blender property, we have (see
(P4') which guaranties (F3) for $f$) that the unstable manifold of
the points $p_i(f)$ is contained in $Q_f$.

Assume that $x \in \partial h_f^{-1}(\{h_f(x)\})$ and consider a
small neighborhood $V$ of $x$. Consider a disk $D$ in
$W^{cs}_{loc}(x)$, since $x$ is a boundary point, we get that
$h_f(D)$ is a compact connected set containing at least two points
in the stable set of $h_f(x)$ for $g$, so by iterating backwards,
and using (F2) (guaranteed for $f$ by (P4')) we get that there is
a backward iterate of $D$ which intersects $W^u_L(p_i(f)) \en Q_f$
for some (i) which concludes.

Now, Theorem B follows with the same argument as for Theorem
\ref{TeoParteTopologica}, using Proposition
\ref{ProposicionMecanismo} and the fact that $r_1(f)$ is contained
in $Q_f$ (because of property (P1')).

\lqqd

\appendix

\section{A new example of robustly transitive attractor}

In this appendix we use a very similar construction to the one
done in Theorem B in order to obtain a robustly transitive
attractor $\Lambda$ in a $3$-dimensional manifold which verify the
following properties:

\begin{itemize}
\item[-] It admits a dominated splitting of the form $TM= E^{cs}
\oplus E^u$ and does not admit any sub-dominated splitting.
\item[-] It appears in any isotopy class of diffeomorphisms of a
$3$-dimensional manifold.
\end{itemize}

We remark that since attractors are invariant under forward
iterates the existence of partially hyperbolic attractors with
splitting of the form $TM=E^{cs} \oplus E^u$ and $TM = E^{s}
\oplus E^{cu}$ is not symmetric. Indeed, it is not hard to show
that the examples of Bonatti and Diaz in \cite{BDPer} can be
modified in order to get robustly transitive attractors with
splitting of the form $TM = E^s \oplus E^{cu}$. Our examples are
by no means completely different, but require a further argument.

Also, it is to be noticed that the example of Carvalho \cite{Car}
also has the same splitting, but since it is a Derived from Anosov
it is not possible to include it in any $3$-manifold.

To sketch the construction of such examples we start with the same
diffeomorphism $F_1$ constructed above.

Let us consider a periodic point $r_1$ of stable index $2$. We can
assume it is fixed (modulo considering an iterate of $F_1$).
Consider $\delta>0$ small enough such that $B(r_1, 6\delta)$ is
disjoint from:

\begin{itemize}
\item[-] the periodic points $p_1, \ldots, p_N, q$ defined above,
\item[-] the blender $K$, \item[-] $\overline{(W^u_L(p_1) \cup
\ldots \cup W^u_L(p_N))}$ and \item[-] from compact connected
pieces of $W^s(p_i)$ ($1\leq i\leq N$) intersecting the activating
region of the blender $K$ and containing $p_i$.
\end{itemize}

\subsubsection{} We shall construct a diffeomorphism $F_2$ modifying $F_1$ such that:

\begin{itemize}
\item[-] $F_2$ coincides with $F_1$ outside $B(r_1,\delta)$.
\item[-] $F_2$ preserves the center-stable lamination of $F_1$.
\item[-] $DF_2$ preserves narrow cones $\cE^u$ and $\cE^{cs}$
around the unstable direction $E^u$ and the center stable
direction $E^{s} \oplus E^c$ of $F_1$ respectively. Also, vectors
in $\cE^{u}$ are expanded uniformly by $DF_2$ while every plane
contained in $\cE^{cs}$ verifies that the volume\footnote{This
means with respect to the Riemannian metric which allows to define
a notion of 2-dimensional volume in each plane.} is contracted by
$DF_2$. \item[-] The point $r_1$ remains fixed for $F_2$ but now
has complex eigenvalues in $r_1$.
\end{itemize}

\begin{prop} There exists an open $C^1$-neighborhood $\cV$ of $F_2$ such that for every $f\in \cV$ one has that $f$ has a transitive attractor in $U$.
\end{prop}

\esbozo Notice that one can choose $\cV$ such that for every $f\in
\cV$ one preserves a center-stable foliation close to the original
one. Also, one can assume that properties (F2) and (F3) still hold
for the continuations $p_i(f)$ and $q(f)$ since $F_2$ coincides
with $F_1$ outside $B(r_1,\delta)$ and these are $C^1$-robust
properties.

Also, we demand that for every $f\in \cV$, the derivative of $f$
preserves the cones $\cE^u$ and $\cE^{cs}$, contracts volume in
$E^{cs}\en \cE^{cs}$ (the plane tangent to the center-stable
foliation) and expands vectors in $E^u \en \cE^u$.

Consider now a center stable disk $D$ and an unstable curve
$\gamma$ which intersect the maximal invariant set

$$Q_f= \bigcap_{n>0} f^n(U) . $$

Since by future iterations $\gamma$ will intersect the stable
manifold of $q(f)$ (property (F3)) we obtain that by the
$\lambda-$lemma it will accumulate the unstable manifold of
$q(f)$. Since the unstable manifold of $q(f)$ is dense in the
union of the unstable manifolds $W^u(p_1(f))\cup \ldots \cup
W^u(p_N(f))$ we obtain that the union of the future iterates of
$\gamma$ will also be dense there.

Now, iterating backwards the disk $D$ we obtain, using that
$Df^{-1}$ expands volume in the center-stable direction that the
diameter of the disk grows exponentially with these iterates.

Condition (F2) will now imply that eventually the backward
iterates of $D$ will intersect the future iterates of $\gamma$.
This implies transitivity. \lqqd


\begin{thebibliography}{2}



\bibitem[A]{A} A. Araujo,  \emph{Exist\^encia de atratores hiperb\'olicos para difeomorfismos de superf\'icie}, Thesis IMPA (1987).



\bibitem[B]{B} C. Bonatti, Towards a global view of dynamical systems, for the $C^1$-topology, \emph{Ergodic Theory and Dynamical Systems.}{\bf 31} n.4 (2011) p. 959-993.

\bibitem[BC]{BC} C. Bonatti and S. Crovisier, R\'ecurrence et G\'en\'ericit\'e,
 \emph{Inventiones Math.}  {\bf 158}  (2004),   33--104.

\bibitem[BCDG]{BCDG} C. Bonatti, S. Crovisier, L. D\'iaz and N. Gourmelon, Internal perturbations of homoclinic classes:non-domination, cycles, and self-replication, \emph{Preprint} arXiv 1011.2935 (2010). To appear in \emph{Ergodic theory and dynamical systems}.

\bibitem[BD1]{BDPer} C. Bonatti and L. D\'iaz, Persistent nonhyperbolic transitive diffeomorphisms, \emph{Ann. Math.} {\bf 143} (1995), 357-396.



\bibitem[BD2]{BD1} C. Bonatti and L. D\'iaz, On maximal generic sets of generic diffeomorphisms, \emph{Publications Math. de l'IHES} {\bf 96} (2003) 171-197.

 \bibitem[BD3]{BD} C. Bonatti and L. D\'iaz, Abundance of $C^1$-robust homoclinic tangencies, \emph{Preprint} arXiv:0909.4062 (2009) to appear in \emph{Transactions of the Amer. Math. Soc.}.

\bibitem[BDP]{BDP} C.Bonatti, L.Diaz and E. Pujals, A $C^1$ generic dichotomy for diffeomorphisms: weak forms of hyperbolicity or infinitely many sinks or sources, \emph{Annals of Math} {\bf 158} (2003), 355--418.

\bibitem[BDV]{BDV} C. Bonatti, L. Diaz and M.Viana, \emph{Dynamics Beyond Uniform Hyperbolicity. A global geometric and probabilistic perspective}, Encyclopaedia of Mathematical Sciences {\bf 102}. Mathematical Physics III. Springer-Verlag (2005).


\bibitem[BLY]{BLY} C.Bonatti, M.Li and D.Yang, On the existence of attractors, \emph{Preprint} arXiv:0904.4393 (2009). To appear in \emph{Transaction of the AMS}.


\bibitem[BV]{BV} C. Bonatti and M.Viana, SRB measures for partially hyperbolic diffeomorphisms whose central
direction is mostly contracting. \emph{Israel J. of Math} {\bf 115} (2000), 157--193.

\bibitem[BF]{BF} J. Buzzi and T. Fisher, Entropic stability beyond partial hyperbolicity, \emph{Preprint} arXiv:1103:2707 (2011).

\bibitem[BFSV]{BFSV} J. Buzzi, T. Fisher, M. Sambarino and C. Vasquez, Maximal entropy measures for certain partially hyperbolic derived from Anosov systems, \emph{Ergodic theory and dynamical systems} (2011).

\bibitem[Car]{Car} M. Carvalho, Sinai-Ruelle-Bowen measures for N-dimensional derived from Anosov diffeomorphisms, \emph{Ergodic Theory and Dynamical Systems}{\bf 13} (1993) 21--44.

\bibitem[C1]{C} S. Crovisier, Partial hyperbolicity far from homoclinic bifurcations, \emph{Advances in Math.} {\bf 226} (2011), 673-726.


\bibitem[C2]{Crov} S. Crovisier, Perturbation de la dynamique de difféomorphismes en topologie $C^1$, \emph{Preprint} arXiv:0912.2896 (2009). To appear in \emph{Asterisque}.



\bibitem[D]{Daverman} R. Daverman, \emph{Decompositions of Manifolds} Academic Press 1986


\bibitem[F]{F} A. Fathi, Expansiveness hyperbolicity and Haussdorf dimension, \emph{Comm. Mathematical Physics} {\bf 126} (1989) 249-262.

\bibitem[HPS]{HPS} M. Hirsch, C. Pugh and M. Shub, Invariant Manifolds,  \emph{Springer Lecture Notes in Math.}, {\bf 583} (1977).




\bibitem[Mi]{Milnor} J. Milnor, On the concept of attractor, \emph{Comm. Math. Physics} {\bf 99} (1985) no.2 177--195.


\bibitem[PV]{PV} J. Palis and M. Viana, High dimensional diffeomorphisms displaying infinitely many periodic attractors, \emph{Annals of Math} {\bf 140} (1994) 207--250.

\bibitem[Pot]{Potrie} R. Potrie, A proof of the existence of attractors in dimension 2, \emph{Unpublished note, not intended for publication} Available at http://www.cmat.edu.uy/$\sim$rpotrie/documentos/pdfs/dimensiondos.pdf .


\bibitem[PS]{PS} E. Pujals and M. Sambarino, Homoclinic tangencies and hyperbolicity for surface diffeomorphisms, \emph{Annals of Math} {\bf 151} (2000) 961--1023.


\bibitem[R]{Rob} C. Robinson, \emph{Dynamical Systems, Stability, Symbolic dynamics, and Chaos} CRC Press (1994).

\bibitem[RHRHTU]{RHRHUT} F. Rodriguez Hertz, M.A. Rodriguez Hertz, A. Tahzibi and R. Ures, Maximizing measures for partially hyperbolic systems with compact center leaves, \emph{Preprint} arXiv 1010.3372 (2010). To appear in \emph{Ergodic theory and dynamical systems}.

\bibitem[VY]{VY} M. Viana and J. Yang, Physical measures and absolute continuity for one-dimensional center direction, \emph{Preprint} IMPA A683 (2010).

\bibitem[W]{W} P. Walters, Anosov diffeomorphisms are topologically stable, \emph{Topology} {\bf 9} (1970) 71--78.

\end{thebibliography}
\end{document}